\documentclass[12pt]{amsart}
 \usepackage{amscd,amsmath,amsthm,amssymb}
 \usepackage{tikz}
  \usepackage{lipsum}
 \usepackage {xcolor}
 \def\NZQ{\mathbb}               

 \def\ZZ{{\NZQ Z}}

 \def\G{{\mathcal G}}

 \def\opn#1#2{\def#1{\operatorname{#2}}} 
 \opn\chara{char} \opn\length{\ell} \opn\pd{pd} \opn\rk{rk}
 \opn\projdim{proj\,dim} \opn\injdim{inj\,dim} \opn\rank{rank}
 \opn\depth{depth} \opn\grade{grade} \opn\height{height}
 \opn\embdim{emb\,dim} \opn\codim{codim}
 
 \opn\Tr{Tr} \opn\bigrank{big\,rank}
 \opn\superheight{superheight}\opn\lcm{lcm}
 \opn\trdeg{tr\,deg}
 \opn\reg{reg} \opn\lreg{lreg} \opn\ini{in} \opn\lpd{lpd}
 \opn\size{size} \opn\sdepth{sdepth}
 \opn\link{link}\opn\fdepth{fdepth}\opn\lex{lex}
 \opn\div{div} \opn\Div{Div} \opn\cl{cl} \opn\Cl{Cl}
\opn\ind{ind}
\opn\match{match}
\opn\deg{deg}
 \opn\Spec{Spec} \opn\Supp{Supp} \opn\supp{supp} \opn\Sing{Sing}
 \opn\Ass{Ass} \opn\Min{Min}\opn\Mon{Mon}
 \opn\Ann{Ann} \opn\Rad{Rad} \opn\Soc{Soc}
 \opn\Im{Im} \opn\Ker{Ker} \opn\Coker{Coker} \opn\Am{Am}
 \opn\Hom{Hom} \opn\Tor{Tor} \opn\Ext{Ext} \opn\End{End}
 \opn\Aut{Aut} \opn\id{id}
 
 \opn\nat{nat}
 \opn\pff{pf}
 \opn\Pf{Pf} \opn\GL{GL} \opn\SL{SL} \opn\mod{mod} \opn\ord{ord}
 \opn\Gin{Gin} \opn\Hilb{Hilb}\opn\sort{sort}
 \opn\Tot{Tot}
 \opn\deg{deg}
\opn\ind{ind}
\opn\match{match}
 \opn\aff{aff} \opn
 \con{conv} \opn\relint{relint} \opn\st{st}
 \opn\lk{lk} \opn\cn{cn} \opn\core{core} \opn\vol{vol}
 \opn\link{link} \opn\star{star}\opn\lex{lex}\opn\set{set}
 \opn\gr{gr}
 
 \def\pot#1#2{#1[\kern-0.28ex[#2]\kern-0.28ex]}

 \opn\dirlim{\underrightarrow{\lim}}
 \opn\inivlim{\underleftarrow{\lim}}
 \let\union=\cup
 \let\sect=\cap

 \let\Dirsum=\bigoplus
 
 \let\to=\rightarrow
 
 \def\Implies{\ifmmode\Longrightarrow \else
         \unskip${}\Longrightarrow{}$\ignorespaces\fi}
 \def\implies{\ifmmode\Rightarrow \else
         \unskip${}\Rightarrow{}$\ignorespaces\fi}
 \def\iff{\ifmmode\Longleftrightarrow \else
         \unskip${}\Longleftrightarrow{}$\ignorespaces\fi}

 \let\:=\colon
 \newtheorem{Theorem}{Theorem}[section]
 
 \newtheorem{Corollary}[Theorem]{Corollary}
 
 \newtheorem{Remark}[Theorem]{Remark}
 
 \newtheorem{Example}[Theorem]{Example}
 
 \newtheorem{Definition}[Theorem]{Definition}

 \let\epsilon\varepsilon
 \let\kappa=\varkappa
 \def\qed{\ifhmode\textqed\fi
       \ifmmode\ifinner\quad\qedsymbol\else\dispqed\fi\fi}
 \def\textqed{\unskip\nobreak\penalty50
        \hskip2em\hbox{}\nobreak\hfil\qedsymbol
        \parfillskip=0pt \finalhyphendemerits=0}
 \def\dispqed{\rlap{\qquad\qedsymbol}}
 
 \opn\dis{dis}
 \def\pnt{{\raise0.5mm\hbox{\large\bf.}}}
 
 \opn\Lex{Lex}

\begin{document}

\title{Edge and vertex cover ideals of planar graphs $St_r$}

\author {Maurizio Imbesi, Monica La\,Barbiera and Roya Moghimipor}

\address{Maurizio Imbesi,  Department of Mathematical and Computer Sciences, Physical and Earth Sciences\\
University of Messina, Viale F. Stagno d'Alcontres 31, I-98166 Messina, Italy} \email{maurizio.imbesi@unime.it}

\address{Monica La\,Barbiera, Department of Electrical, Electronic and Computer Engineering\\
 University of Catania, Viale A. Doria 6, I-95125 Catania, Italy} \email{monica.labarbiera@unict.it}

\address{Roya Moghimipor, Department of Mathematics, Safadasht Branch\\
Islamic Azad University, Tehran, Iran}  \email{roya\_moghimipour@yahoo.com}

\begin{abstract}
The interplay between topology, graph theory and commutative algebra represents one of the most prolific research directions in modern algebraic and geometric combinatorics. In this paper, we conduct an extensive structural, homological and algebraic investigation of the edge ideal $I(St_r)$ and the vertex cover ideal $I_c(St_r)$ associated with the family of graphs $St_r\,, \,r > 1$ the number of their plane regions, $2r+1$ the number of their vertices.\\
In the first part, we focus on $I_c(St_r)$ and its powers $I_c(St_r)^k, \,k\geq 1$ arbitrary integer. By establishing the explicit irredundant primary decomposition of $I_c(St_r)^k$, we show that $\text{height}(I_c(St_r)^k) = 2$ holds for every $k$, which systematically yields the result that $I_c(St_r)^k$ is an unmixed ideal across all powers $k$. Furthermore, by constructing inductive minimal $\mathbb{Z}^{2r+1}$-graded free resolutions for the quotient ring $R/I_c(St_r)$, we derive closed-form piecewise linear formulas for two fundamental homological invariants: the projective dimension $\text{proj\,dim}(R/I_c(St_r))$ and the Castelnuovo--Mumford regularity $\text{reg}(R/I_c(St_r))$, both classified explicitly according to $r \pmod 3$.\\
In the second part, we turn our attention to syzygy modules and free resolutions of $I(St_r)$. We calculate the exact number of minimal quadratic generators, and the number of linear syzygies, by connecting $I(St_r)$ to the structure of the edge graph $\mathcal{L}(St_r)$ of $St_r$. Finally, we address a well-known conjecture by Eliahou and Villarreal which concerns the second syzygy module and the third graded Betti number in degree 4. We validate and prove this conjecture for several fundamental classes of graphs, including graphs $St_r$.
\end{abstract}

\thanks{}

\maketitle

\section*{Introduction}
In modern commutative algebra, the study of monomial ideals associated with combinatorial structures has emerged as an essential tool for understanding the bridge between discrete mathematics and homological algebra. The foundation of this field relies on Stanley--Reisner theory, which establishes a natural dictionary between simplicial complexes and squarefree monomial ideals \cite{HerzogHibi2010, Villarreal2015}. Given a finite simple graph $G = (V(G), E(G))$ where $V(G) = \{v_1, \dots, v_n\}$ denotes its vertex set and $E(G)$ its edge set, one can associate squarefree monomial ideals in a polynomial ring $R = K[x_1, \dots, x_n]$ over a field $K$, where each variable $x_i$ uniquely corresponds to a vertex $v_i$ of $G$.

The key algebraic object attached to $G$ is its \emph{edge ideal} $I(G) = (x_i x_j \mid \{v_i, v_j\} \in E(G))$, introduced by Villarreal \cite{Villarreal2015}. The edge ideal directly encapsulates the adjacency relations of the graph. The dual to edge ideal is the \emph{vertex cover ideal} $I_c(G)$ (just studied via Alexander duality), which is generated by squarefree monomials corresponding to minimal vertex covers of $G$ \cite{HerzogHibi2010, ImbesiLaBarbiera2015}. A subset $\mathcal{C}(G) \subseteq V(G)$ is a vertex cover if every edge in $E(G)$ has at least one endpoint in $\mathcal{C}(G)$. The minimal prime ideals of $I(G)$ correspond bijectively to the minimal vertex covers of $G$, meaning that $I_c(G)$ encodes the primary decomposition of $I(G)$ and provides critical information regarding its height, primary components, and associated prime ideals \cite{ImbesiLaBarbiera2014, Villarreal2007}.

A graph $G$ is called \emph{unmixed} if all its minimal vertex covers have the same cardinality, or equivalently, if $R/I(G)$ is an unmixed quotient ring \cite{Villarreal2007}. A major topic of current interest in combinatorial commutative algebra is the analysis of powers of monomial ideals; specifically, whether algebraic properties such as unmixedness, Cohen--Macaulayness, or linear quotients are preserved under taking ordinary powers $I_c(G)^k$ or symbolic powers $I_c(G)^{(k)}$ for $k \ge 1$ \cite{ImbesiLaBarbiera2014, ImbesiLaBarbiera2015}.

Planar graphs form a classical family in topological graph theory. By Kuratowski's theorem, a graph is planar if and only if it contains no subgraph isomorphic to or homeomorphic to the complete graph $\mathcal{K}_5$ or the complete bipartite graph $\mathcal{K}_{3,3}$ \cite{Kuratowski1930}. Due to their spatial constraints, planar graphs exhibit rigid upper bounds on their number of edges relative to their vertex count \cite{BondyMurty2008, Harary1972}. In 1996, Doering and Gunston introduced a specific infinite family of bipartite planar graphs, denoted by $St_r$ for $r > 1$, which exhibit rich geometric and algebraic symmetries \cite{DoeringGunston1996}. The graph $St_r$ consists of $r$ planar regions on the vertex set $V(St_r) = \{v_1, \dots, v_{2r+1}\}$, partitioned into disjoint subsets $V_1 = \{x_1, \dots, x_{r+1}\}$ and $V_2 = \{y_1, \dots, y_r\}$, with an edge set defined by:
\begin{eqnarray*}
&E(St_{r})=&\{\{v_1,v_i\} \mid  2 \leq i \leq r + 1\} \cup \{\{v_i,v_{i+r}\} \mid  2 \leq i \leq r + 1\}
\\
&&
\cup \{\{v_i,v_{i+r-1}\} \mid  3\leq i \leq r+1\}\cup \{v_2,v_{2r+1}\}.
\end{eqnarray*}
When $r = 2$, $St_2$ coincides with the complete bipartite graph $\mathcal{K}_{3,2}$, whereas for $r > 2$, $St_r$ is a sparse bipartite planar graph \cite{DoeringGunston1996}.

The primary goal of this paper is to conduct a detailed investigation into both the cover ideal $I_c(St_r)$ and the edge ideal $I(St_r)$, bridging their combinatorial graph properties with homological invariants such as Betti numbers, projective dimension, and Castelnuovo--Mumford regularity \cite{Jacques2004, EliahouVillarreal1999}.

The manuscript is organized into two main parts.

In Section 1, we examine the algebraic behavior of the vertex cover
ideal $I_c(St_r)$ and its ordinary powers $I_c(St_r)^k$. First, we
compute the irredundant primary decomposition of $I_c(St_r)^k$ for
any integer $k \ge 1$. We prove that $\text{height}(I_c(St_r)^k) =
2$ for all $k \ge 1$, which immediately implies that $I_c(St_r)^k$
is unmixed for every $k \ge 1$. Next, by explicitly computing the
minimal $\mathbb{Z}^{2r+1}$-graded free resolutions for
$R/I_c(St_r)$ via an inductive argument on the number of regions
$r$, we derive exact closed formulas for both the projective
dimension $\text{proj\,dim}(R/I_c(St_r))$ and the
Castelnuovo--Mumford regularity $\text{reg}(R/I_c(St_r))$. We
establish that these invariants behave as piecewise linear functions
depending on $r \pmod 3$, building upon general homological theory
of monomial ideals \cite{HerzogHibi2010}.

In Section 2, we shift our focus to the syzygy modules and free
resolutions of the edge ideal $I(St_r)$. Using properties of
edge graphs $\mathcal{L}(G)$ of a graph $G$ and formulas from Eliahou and
Villarreal \cite{EliahouVillarreal1999}, we determine the exact
number of minimal generators $q = 3r$ and the number of linear
syzygies $b = \frac{1}{2}r(r+7)$ for $I(St_r)$. Furthermore, we
investigate a prominent conjecture raised in
\cite{EliahouVillarreal1999} concerning the third graded Betti
number in degree 4, $b_{3,4}(R/I(G)) = d$, which represents the
number of linear generators in the second syzygy module of $I(G)$.
We provide comprehensive proofs of this conjecture for complete graphs, complete bipartite graphs, cycle graphs, and for the class of bipartite planar graphs $St_r$, proving that in all
these cases
\[
d = \sum_{i=1}^{N} \binom{\text{deg}(v_i)}{3} + W_4,
\]
where $N$ is the number of vertices and $W_4$ the number of chordless 4-cycles (\cite{Jacques2004}).

\bigskip

\section{Invariants of vertex cover ideals of graphs $St_r$}
\label{one}

In this section we illustrate some algebraic aspects of the vertex cover ideals of planar graphs.

Let $R=K[x_1,\dots,x_n]$ be a polynomial ring in $n$ variables over an arbitrary field $K$ such that
$\deg x_i = 1$, for $i = 1,\dots,n$.
For any monomial ideal $I$ of $R$, let $\G(I)$ denote the unique minimal set of monomial generators for $I$.

\begin{Definition}
{\em
A {\em vertex cover} of $I \subset R$ is a subset $\mathcal{C}$ of $\{x_1,\dots, x_n\}$ such that each $u \in \G(I)$ is divided by
some $x_i \in \mathcal{C}$. Then $\mathcal{C}$ is called {\em minimal} if no proper subset of $\mathcal{C}$ is a vertex cover of $I$.
We denote the minimal cardinality of the vertex covers of $I$ by $\height(I)$.
}
\end{Definition}

Let us now give information on graphs and their ideals.
Let $G$ be a graph on the vertex set $V(G)=\{v_1,\dots,v_n\}$ and edge set $E(G)$,
that consists of pairs $\{v_i,v_j\}$ said edges, for some $v_i,v_j \in V(G)$.

A graph $G$ on $n$ vertices $v_1,\dots, v_n$ is {\em complete} if there exists an edge for all
pairs $\{v_i,v_j\}$ of vertices of $G$. It is denoted by $\mathcal{K}_{n}$.
$G$ is {\em bipartite} if its vertex set $V(G)$ can be partitioned into disjoint subsets $V_1=\{x_1,\dots,x_n\}$
and $V_2=\{y_1,\dots,y_m\}$, and any edge joins a vertex of $V_1$ to a vertex of $V_2$.
A bipartite graph $G$ is {\em complete} if all the vertices of $V_1$ are joined to all the
vertices of $V_2$. It is denoted by $\mathcal{K}_{n,m}$.

Let $G$ be a graph on vertices $v_1,\dots,v_n$ and $R=K[x_1,\dots,x_n]$ be a
polynomial ring over a field $K$, where each variable $x_i$ corresponds to vertex $v_i$.

\begin{Definition}
{\em
Let $G$ be a graph. The {\em edge ideal} $I(G)$ of $G$ is the squarefree monomial ideal of
$R$ generated by squarefree monomials of degree two, $x_ix_j$, on the $x_1,\dots,x_n$ variables,
such that $\{v_i,v_j\}\in E(G)$ for $1 \leq i, j \leq n$:
\[
I(G)=(x_ix_j \mid \{v_i,v_j\}\in E(G)).
\]
}
\end{Definition}

\begin{Definition}
{\em
A graph $G$ is {\em planar} if it has an embedding in the plane
such that each pair of edges is intersected in common vertices only.
A planar graph is divided by its edges in plane regions.
}
\end{Definition}

\begin{Remark}
{\em
The complete graphs $\mathcal{K}_5$ and $\mathcal{K}_{3,3}$ are the minimal not planar
graphs. In fact it is not possible to represent these graphs in the plane so
that edges are not intersected in common vertices only.
}
\end{Remark}

\begin{Theorem}
\label{planar}
A graph is planar if and only if it has no subgraphs containing $\mathcal{K}_{5}$ and $\mathcal{K}_{3,3}$.
\end{Theorem}
\begin{proof}
Take into consideration that one in \cite{Harary1972}, Theorem 11.13\,.
\end{proof}

\begin{Remark}
\label{condition}
{\em
(\cite{BondyMurty2008, Harary1972}) Let $G$ be a planar graph, and let $N$ be the number of its vertices, $w$ be the
number of its edges. The following conditions hold:
\begin{itemize}
\item[(1)] If $N = 3$, then $w=3N-6$;

\item[(2)] If $N > 3$ and there are not cycles of length $3$, then $w=2N-4$.
\end{itemize}
}
\end{Remark}

\begin{Example}
{\em
The graph $\mathcal{K}_{3,3}$ has $N=6$ vertices and no cycle of length $3$, but $w=9$.
Hence $\mathcal{K}_{3,3}$ does not satisfy the condition (2) of Remark \ref{condition}. Therefore,
$\mathcal{K}_{3,3}$ is not a planar graph.
}
\end{Example}

We now examine the class of planar graphs $St_{r}$ introduced in \cite{DoeringGunston1996}.

\noindent Let $St_{r}$ be the planar graph with $r > 1$ regions on vertex set
$V(St_{r})=\{v_1,\dots,v_{2r+1}\}$ and edge set
\begin{eqnarray*}
&E(St_{r})=&\{\{v_1,v_i\} \mid  2 \leq i \leq r + 1\} \cup \{\{v_i,v_{i+r}\} \mid  2 \leq i \leq r + 1\}
\\
&&
\cup \{\{v_i,v_{i+r-1}\} \mid  3\leq i \leq r+1\}\cup \{v_2,v_{2r+1}\}.
\end{eqnarray*}

\begin{Example}
\label{St2}
{\em
$St_{2}$ is the bipartite planar graph on vertex set $V(St_{2})=\{v_1,\dots,v_5\}$ and
edge set
\[
E(St_{2})=\{\{v_1,v_2\},\{v_1,v_3\},\{v_2,v_4\},\{v_3,v_5\},\{v_3,v_4\},\{v_2,v_5\}\}.
\]
Hence, $V(St_{2})$ can be partitioned into disjoint subsets
$
V_1=\{x_1,x_2,x_3\}$ and $V_2=\{y_1,y_2\}
$
such that $x_1=v_1, x_2=v_4, x_3=v_5, y_1=v_2, y_2=v_3$.
Therefore,
\[
E(St_{2})=\{\{x_1,y_1\}, \{x_1,y_2\}, \{x_2,y_1\}, \{x_3,y_2\}, \{x_2,y_2\}, \{x_3,y_1\}\},
\]
and hence $St_{2}$ is a complete bipartite planar graph.\\

\begin{center}
\begin{minipage}{.4\textwidth}
\usetikzlibrary {positioning}
\definecolor {processblue}{cmyk}{0.96,0,0,0}
\begin{tikzpicture}
      \tikzset{enclosed/.style={draw, circle, inner sep=0pt, minimum size=.15cm, fill=black, text=blue , minimum width =0.1 cm,node distance =0.1 cm,top color =blue , bottom color = processblue!20 ,
draw,processblue , text=blue , minimum width =0.1 cm }}
\node[enclosed, label={left, yshift=.2cm: $x_{1}$}] (v_{1}) at (2,6) {};
\node[enclosed, label={left, yshift=.2cm: $x_{2}$}] (v_{4}) at (4,6) {};
\node[enclosed, label={left, yshift=.2cm: $x_{3}$}] (v_{5}) at (6,6) {};

\node[enclosed, label={left, yshift=.1cm: $y_{1}$}] (v_{2}) at (3,4) {};
\node[enclosed, label={left, yshift=.01cm: $y_{2}$}] (v_{3}) at (5,4) {};

\draw (v_{1}) -- (v_{2}) node[midway, right] (edge1){};
\draw (v_{1}) -- (v_{3}) node[midway, below] (edge2){};
\draw (v_{2}) -- (v_{4}) node[midway, below] (edge3){};
\draw (v_{4}) -- (v_{3}) node[midway, below] (edge4){};
\draw (v_{5}) -- (v_{3}) node[midway, below] (edge5){};
\draw (v_{4}) -- (v_{3}) node[midway, below] (edge6){};
\draw (v_{5}) -- (v_{2}) node[midway, below] (edge7){};

\end{tikzpicture}
\end{minipage}
\end{center}

\medskip

\begin{center}
\begin{minipage}{.4\textwidth}
\usetikzlibrary {positioning}
\definecolor {processblue}{cmyk}{0.96,0,0,0}
\begin{tikzpicture}
      \tikzset{enclosed/.style={draw, circle, inner sep=0pt, minimum size=.15cm, fill=black, text=blue , minimum width =0.1 cm,node distance =0.1 cm,top color =blue , bottom color = processblue!20 ,
draw,processblue , text=blue , minimum width =0.1 cm }}
\node[enclosed, label={left, yshift=.2cm: $v_{5}$}] (v_{5}) at (2,5) {};
\node[enclosed, label={left, yshift=-.05cm: $v_{3}$}] (v_{3}) at (4,4) {};
\node[enclosed, label={left, yshift=.3cm: $v_{1}$}] (v_{1}) at (3,2.5) {};
\node[enclosed, label={left, yshift=.2cm: $v_{2}$}] (v_{2}) at (1,3) {};
\node[enclosed, label={left, yshift=-.08cm: $v_{4}$}] (v_{4}) at (4,1) {};

\draw (v_{1}) -- (v_{3}) node[midway, right] (edge1){};
\draw (v_{1}) -- (v_{2}) node[midway, below] (edge2){};
\draw (v_{2}) -- (v_{5}) node[midway, below] (edge3){};
\draw (v_{3}) -- (v_{5}) node[midway, below] (edge4){};
\draw (v_{2}) -- (v_{4}) node[midway, below] (edge5){};
\draw (v_{3}) -- (v_{4}) node[midway, below] (edge6){};

\end{tikzpicture}
\end{minipage}
\end{center}
\smallskip
The two pictures represent the same graph $St_{2}$.
}
\end{Example}

Now we want to study algebraic and homological properties of powers of vertex cover ideals of $St_{r}$.

\begin{Definition}
{\em
Let $G$ be a graph on vertex set $V(G)$. A subset $\mathcal{C}(G)\subset V(G)$ is
called a \textit{minimal vertex cover} for $G$ if:\\
(1) each edge of $G$ is incident with one vertex in $\mathcal{C}(G)$; \\
(2) there is no proper subset  of $\mathcal{C}(G)$ with this property.
}
\end{Definition}

If $\mathcal{C}(G)$ satisfies condition (1) only, then $\mathcal{C}(G)$ is called a \textit{vertex cover} of $G$ and it is said to cover all the edges of $G$.\\
\noindent Therefore, the smallest number of vertices in any  minimal vertex cover of $G$ is called \textit{vertex covering number}.\\[-3mm]

There exists a one to one correspondence between minimal vertex covers of G and minimal prime
ideals of the edge ideal $I(G)$. Hence, $\mathcal{P}$ is a minimal prime ideal of $I(G)$ if and only if $\mathcal{P}=(\mathcal{C}(G))$, for some minimal vertex cover $\mathcal{C}(G)$ of $G$.
Therefore, $I(G)$ has primary decomposition $(\mathcal{C}_1(G)),\dots,(\mathcal{C}_q(G))$,
where $\mathcal{C}_1(G),\dots,\mathcal{C}_q(G)$ are the minimal vertex covers of $G$.

\begin{Definition}
{\em
Let $G$ be a graph. The {\em vertex cover ideal} of $G$, denoted by $I_c(G)$, is the monomial ideal of $R$ generated by all the monomials $x_{i1}\dots x_{iw}$ such that $\{x_{i1},\dots,x_{iw}\}$ is an associated prime ideal of $I(G)$.
}
\end{Definition}

\smallskip

In the following, we compute the height of powers of the vertex cover ideal $I_{c}(St_{r})$.

\begin{Theorem}
\label{height}
Let $R=K[x_1,\dots,x_n]$ and $I_{c}(St_{r})\subset R$. Then \vspace{2mm}\par
\hspace{3cm} $\height(I_{c}(St_{r})^{k})=2, \;\;$ for every $k\geq 1$.
\end{Theorem}

\begin{proof}
We suppose that $St_{r}$ be the planar graph with $r > 1$ regions on vertex set
$V(St_{r})=\{v_1,\dots,v_{2r+1}\}$, and edge set
\begin{eqnarray*}
&E(St_{r})=&\{\{v_1,v_i\} \mid  2 \leq i \leq r + 1\} \cup \{\{v_i,v_{i+r}\} \mid  2 \leq i \leq r + 1\}
\\
&&
\cup \{\{v_i,v_{i+r-1}\} \mid  3\leq i \leq r+1\}\cup \{v_2,v_{2r+1}\}.
\end{eqnarray*}

The vertex set of $St_{r}$ can be partitioned
into disjoint subsets $V_1=\{x_1,\dots,x_n\}$ and $V_2=\{y_1,\dots,y_m\}$,
where $n=r+1$ and $m=r$, $m+n=2r+1$, and any edge joins a vertex of
$V_1$ with a vertex of $V_2$. Then the edge set can be written:
\begin{eqnarray*}
&E(St_{r})=&\{\{x_1,y_i\} \mid 1 \leq i \leq m\} \cup \{\{x_i,y_{i-1}\} \mid 2 \leq i \leq n \}
\\
&&
\cup
\{\{x_i,y_i\} \mid 2 \leq i \leq m\} \cup \{x_n,y_1\}.
\end{eqnarray*}
It follows that $St_{r}$ is a bipartite planar graph, and it is a complete bipartite planar graph in the case $r=2$ only.
Therefore, the irredundant primary decomposition of $I_{c}(St_{r})^{k}$ is the following:
\begin{eqnarray*}
&I_{c}(St_{r})^{k}=&\left(\bigcap_{1 \leq i \leq m}(x_1,y_i)^{k}\right)\sect \left(\bigcap_{2 \leq i \leq n}(x_i,y_{i-1})^{k}\right)\sect
\\
&&
\left(\bigcap_{2 \leq i \leq m}(x_i,y_i)^{k}\right) \sect (x_n,y_1)^{k},
\end{eqnarray*}
for every $k\geq 1$. Hence, $\height(I_{c}(St_{r})^{k})=2$ for every $k\geq 1$, as desired.
\end{proof}

We recall that a graph $G$ is {\em unmixed} if and only if all the minimal vertex covers of $G$ have the same cardinality. From Theorem \ref{height} the following result descends:

\begin{Corollary}
\label{unmixed}
Let $St_{r}$ be the bipartite planar graph with $r > 1$ regions, and let $I_{c}(St_{r})$ be its vertex cover ideal.
Then $I_{c}(St_{r})^{k}$ is unmixed for every $k\geq 1$.
\end{Corollary}

\begin{proof}
We assume that $St_{r}$ be the bipartite planar graph with $r > 1$ regions.
The vertex set of $St_{r}$ can be partitioned into disjoint subsets $V_{1}$ and $V_{2}$
such that $|V_{1}|+|V_{2}|=2r+1$, and $|V_{i}|$ denotes the number of vertices of $V_{i}$ for $i=1,2$.
Hence, by Theorem \ref{height} we conclude that the height of all minimal prime ideals of $I_{c}(St_{r})^{k}$
are equal to $2$. Then $I_{c}(St_{r})^{k}$ is unmixed, and the assertion follows.
\end{proof}

Moreover, let's give formulas for the projective dimension and regularity of $R/I_{c}(St_{r})$.\\[2mm]
Let $M$ be a finitely generated graded $R$-module. Hence, $M$ has a minimal graded free resolution of the form
\[
 0\to \Dirsum_{g\in \ZZ}S(-g)^{b_{p,g}(M)}\to \cdots \to \Dirsum_{g\in \ZZ}S(-g)^{b_{1,g}(M)}\to \Dirsum_{g\in \ZZ}S(-g)^{b_{0,g}(M)}\to M\to 0,
\]
where $b_{i,g}(M)$ are called the {\em graded Betti numbers} of $M$.\\ 
The {\em projective dimension} of $M$, denoted by $\projdim(M)$, is defined by
\[
\projdim(M)=\max\{i \ | \ b_{i,g}(M)\neq 0 \;\; \text{for some $g$}\}.
\]
The {\em Castelnuovo-Mumford regularity} of $M$, denoted by $\reg(M)$, is defined by
\[
\reg(M)=\max\{g-i \ | \ b_{i,g}(M)\neq 0\}.
\]
For more information we refer the reader to \cite{HerzogHibi2010}.

\begin{Theorem}
\label{projdim}
Let $St_{r}$ be a bipartite planar graph with $r > 1$ regions, and let $I_{c}(St_{r})$ be its vertex cover ideal. Then
$$ \projdim(R/I_{c}(St_{r}))=\left\{ \begin{array}{ll}
\frac{2}{3}r+\frac{2}{3} & \quad \textrm{if \quad  $r=3w-1$},\\[2mm]
\frac{2}{3}r+1 &  \quad \textrm{if \quad $r=3w$},\\[2mm]
\frac{2}{3}r+\frac{4}{3}  & \quad \textrm{if \quad $r=3w+1$},\\
\end{array} \right.
$$
for every $w \geq1$.
\end{Theorem}

\begin{proof}
Let $I_{c}(St_{r})$ be the vertex cover ideal of the bipartite planar graph $St_{r}$ with $r > 1$ regions.
Furthermore, we assume that $w$ is a positive integer.

The vertex set of $St_{r}$ can be partitioned into disjoint subsets $V_{1}$ and $V_{2}$
such that $|V_{1}|+|V_{2}|=2r+1$, and $|V_{i}|$ denotes the number of vertices of $V_{i}$ for $i=1,2$.
\noindent We distinguish three cases.

\medskip

$\mathbf{Case(\.{i}):}$ Let $r=3w-1$. We prove the assertion by induction on $r$.
\begin{itemize}
\item If $r=2$, then
$
V(St_{2})=\{v_1,\dots,v_5\}
$
and edge set
\end{itemize}
\[
E(St_{2})=\{\{v_1,v_2\},\{v_1,v_3\},\{v_2,v_4\},\{v_3,v_5\},\{v_3,v_4\},\{v_2,v_5\}\}.
\]
The vertex cover ideal $I_{c}(St_{2})$ has $2$ monomial squarefree generators and it is:
\[
I_{c}(St_{2})=(x_{2}x_{3},x_{1}x_{4}x_{5})
\]
is the vertex cover ideal of $St_{2}$.
Then the minimal $\ZZ^{5}$-graded free resolution of $R/I_{c}(St_{2})$ is the following:
\[
0 \rightarrow R(-5)\rightarrow R(-2)\oplus R(-3)\rightarrow R\rightarrow R/I_{c}(St_{2})\rightarrow 0
\]
It follows that
$
\projdim(R/I_{c}(St_{2}))=2.
$
\begin{itemize}
\item If $r=5$, hence $St_{5}$ is the bipartite planar graph on the vertex set
\end{itemize}
$
V(St_{5})=\{v_1,\dots,v_{11}\}
$
and edge set
\begin{eqnarray*}
&E(St_{5})=&\{\{v_1,v_i\} \mid  2 \leq i \leq 6\} \cup \{\{v_i,v_{i+5}\} \mid  2 \leq i \leq 6\}
\\
&&
\cup \{\{v_i,v_{i+4}\} \mid  3\leq i \leq 6\}\cup \{v_2,v_{11}\}.
\end{eqnarray*}
Therefore,
\[
0\rightarrow R(-10) \rightarrow R^{15}(-9) \rightarrow R^{30}(-8) \rightarrow
\]
\[
R(-5)\oplus R(-6) \oplus R^{15}(-7)\rightarrow R\rightarrow R/I_{c}(St_{5})\rightarrow 0
\]
is the minimal $\ZZ^{11}$-graded free resolution of $R/I_{c}(St_{5})$.
Hence,
\[
\projdim(R/I_{c}(St_{5}))=\projdim(R/I_{c}(St_{2}))+ \projdim(R/I_{c}(St_{2}))=4.
\]

\begin{itemize}
\item If $r=8$, then $V(St_{8})=\{v_1,\dots,v_{17}\}$ and edge set
\begin{eqnarray*}
&E(St_{8})=&\{\{v_1,v_i\} \mid  2 \leq i \leq 9\} \cup \{\{v_i,v_{i+8}\} \mid  2 \leq i \leq 9\}
\\
&&
\cup \{\{v_i,v_{i+7}\} \mid  3\leq i \leq 9\}\cup \{v_2,v_{17}\}.
\end{eqnarray*}
\end{itemize}
The minimal $\ZZ^{17}$-graded free resolution of $R/I_{c}(St_{8})$ is the following:
\[
0 \rightarrow  R(-17) \rightarrow R^{40}(-15) \rightarrow R^{160}(-14)\rightarrow R^{48}(-12)\oplus R^{240}(-13) \rightarrow
\]
\[
R^{96}(-11)\oplus R^{160}(-12) \rightarrow R(-8) \oplus R(-9)\oplus R^{40}(-10)\oplus R^{40}(-11) \rightarrow R\rightarrow R/I_{c}(St_{8})\rightarrow 0
\]
It follows that
\[
\projdim(R/I_{c}(St_{8}))=\projdim(R/I_{c}(St_{5}))+ \projdim(R/I_{c}(St_{2}))=6.
\]
Therefore, for every $w \geq1$ we have

\setlength{\tabcolsep}{1.5pt}
\renewcommand{\arraystretch}{1}
\begin{center}
\begin{tabular}{|c|c|}
\hline
$r=3w-1$&$\projdim R/I_{c}(St_{r})$\\[1mm]
\hline
2&2\\
\hline
5&$\projdim(R/I_{c}(St_{2}))+ \projdim(R/I_{c}(St_{2})=2+2=4$\\
\hline
8&$\projdim(R/I_{c}(St_{5}))+ \projdim(R/I_{c}(St_{2})=4+2=6$\\
\hline
11&$\projdim(R/I_{c}(St_{8}))+ \projdim(R/I_{c}(St_{2})=6+2=8$\\
\hline
14&$\projdim(R/I_{c}(St_{11}))+ \projdim(R/I_{c}(St_{2})=8+2=10$\\
\hline
17&$\projdim(R/I_{c}(St_{14}))+ \projdim(R/I_{c}(St_{2})=10+2=12$\\
\hline
20&$\projdim(R/I_{c}(St_{17}))+ \projdim(R/I_{c}(St_{2})=12+2=14$\\
\hline
$\vdots$&\vdots\\
\hline
\end{tabular}
\end{center}

\medskip

\noindent Assume that the induction hypothesis holds for $r-3$. It follows that
\begin{eqnarray*}
\projdim(R/I_{c}(St_{3w-1}))&=&\projdim(R/I_{c}(St_{3w-1-3}))+ \projdim(R/I_{c}(St_{2}))\\
&=&2(w-1)+2=2w,
\end{eqnarray*}
for every $w \geq2$.
Hence, if $r=3w-1$,\,
$
\projdim(R/I_{c}(St_{r}))=\frac{2}{3}r+\frac{2}{3}
$
for \mbox{every $w \geq 1$.}

\medskip

$\mathbf{Case(\.{i}\.{i}):}$ Let $r=3w$. By induction on $r$, we prove the assertion.

\begin{itemize}
\item If $r=3$, hence $St_{3}$ is the bipartite planar graph on the vertex set
\end{itemize}
$
V(St_{3})=\{v_1,\dots,v_{7}\}
$
and edge set
\begin{eqnarray*}
&E(St_{3})=&\{\{v_1,v_i\} \mid  2 \leq i \leq 4\} \cup \{\{v_i,v_{i+3}\} \mid  2 \leq i \leq 4\}
\\
&&
\cup \{\{v_i,v_{i+2}\} \mid  3\leq i \leq 4\}\cup \{v_2,v_{7}\}.
\end{eqnarray*}
Therefore, the minimal $\ZZ^{7}$-graded free resolution of $R/I_{c}(St_{3})$ is the following:
\[
0 \rightarrow R^{2}(-7)\rightarrow R^{6}(-6) \rightarrow R(-3)\oplus R(-4) \oplus R^{3}(-5)\rightarrow R\rightarrow R/I_{c}(St_{3})\rightarrow 0
\]
It follows that
$
\projdim(R/I_{c}(St_{3}))=3.
$

\begin{itemize}
\item If $r=6$, then $St_{6}$ is the bipartite planar graph on the vertex set
\end{itemize}
$
V(St_{6})=\{v_1,\dots,v_{13}\}
$
and edge set
\begin{eqnarray*}
&E(St_{6})=&\{\{v_1,v_i\} \mid  2 \leq i \leq 7\} \cup \{\{v_i,v_{i+6}\} \mid  2 \leq i \leq 7\}
\\
&&
\cup \{\{v_i,v_{i+5}\} \mid  3\leq i \leq 7\}\cup \{v_2,v_{13}\}.
\end{eqnarray*}
Therefore, the minimal $\ZZ^{13}$-graded free resolution of $R/I_{c}(St_{6})$ is the following:
\[
0 \rightarrow R^{2}(-13)\rightarrow
R^{12}(-12)\rightarrow R^{24}(-10)\oplus R^{18}(-11) \rightarrow
\]
\[
R^{48}(-9)\oplus R^{12}(-10) \rightarrow R(-6) \oplus R(-7)\oplus R^{24}(-8)\oplus R^{3}(-9) \rightarrow R\rightarrow R/I_{c}(St_{6})\rightarrow 0
\]
Hence,
$
\projdim(R/I_{c}(St_{6}))=\projdim(R/I_{c}(St_{3}))+ (\projdim(R/I_{c}(St_{3}))-1)=5.
$
Therefore, if $r=3w$ we have

\setlength{\tabcolsep}{1.5pt}
\renewcommand{\arraystretch}{1}
\begin{center}
\begin{tabular}{|c|c|}
\hline
$r=3w$&$\projdim(R/I_{c}(St_{r}))$\\[1mm]
\hline
3&3\\
\hline
6&$\projdim(R/I_{c}(St_{3}))+ (\projdim(R/I_{c}(St_{3}))-1)=3+2=5$\\
\hline
9&$\projdim(R/I_{c}(St_{6}))+ (\projdim(R/I_{c}(St_{3}))-1)=5+2=7$\\
\hline
12&$\projdim(R/I_{c}(St_{9}))+ (\projdim(R/I_{c}(St_{3}))-1)=7+2=9$\\
\hline
15&$\projdim(R/I_{c}(St_{12}))+ (\projdim(R/I_{c}(St_{3}))-1)=9+2=11$\\
\hline
18&$\projdim(R/I_{c}(St_{15}))+ (\projdim(R/I_{c}(St_{3}))-1)=11+2=13$\\
\hline
21&$\projdim(R/I_{c}(St_{18}))+ (\projdim(R/I_{c}(St_{3}))-1)=13+2=15$\\
\hline
$\vdots$&\vdots\\
\hline
\end{tabular}
\end{center}

\medskip

\noindent Assume that the induction hypothesis holds for $r-3$. Then
\begin{eqnarray*}
\projdim(R/I_{c}(St_{3w}))&=&\projdim(R/I_{c}(St_{3w-3}))+ \left(\projdim(R/I_{c}(St_{3}))-1\right)\\
&=&\left(2(w-1)+1\right)+(3-1)\\
&=&2w+1.
\end{eqnarray*}
Therefore, if $r=3w$, \,$\projdim(R/I_{c}(St_{r}))=\frac{2}{3}r+1$ for every $w \geq1$.

\medskip

$\mathbf{Case(\.{i}\.{i}\.{i}):}$ Let $r=3w+1$ for every $w \geq1$. Then by induction on $r$,
we prove that $\projdim(R/I_{c}(St_{r}))=2w+2$ for every $w \geq1$.

\begin{itemize}
\item If $r=4$, hence $St_{4}$ is the bipartite planar graph on the vertex set
\end{itemize}
$
V(St_{4})=\{v_1,\dots,v_{9}\}
$
and edge set
\begin{eqnarray*}
&E(St_{4})=&\{\{v_1,v_i\} \mid  2 \leq i \leq 5\} \cup \{\{v_i,v_{i+4}\} \mid  2 \leq i \leq 5\}
\\
&&
\cup \{\{v_i,v_{i+3}\} \mid  3\leq i \leq 5\}\cup \{v_2,v_{9}\}.
\end{eqnarray*}
Therefore,
\[
0 \rightarrow R(-9) \rightarrow R^{8}(-8) \rightarrow R^{16}(-7) \rightarrow R(-4)\oplus R(-5)\oplus R^{8}(-6) \rightarrow R\rightarrow R/I_{c}(St_{4})\rightarrow 0
\]
is the minimal $\ZZ^{9}$-graded free resolution of $R/I_{c}(St_{4})$.
It follows that
\[
\projdim(R/I_{c}(St_{4}))=4.
\]

\begin{itemize}
\item If $r=7$, then $V(St_{7})=\{v_1,\dots,v_{15}\}$ and edge set
\end{itemize}
\begin{eqnarray*}
&E(St_{7})=&\{\{v_1,v_i\} \mid  2 \leq i \leq 8\} \cup \{\{v_i,v_{i+7}\} \mid  2 \leq i \leq 8\}
\\
&&
\cup \{\{v_i,v_{i+6}\} \mid  3\leq i \leq 8\}\cup \{v_2,v_{15}\}.
\end{eqnarray*}
It follows that
\[
0\rightarrow R(-16) \rightarrow R^{36}(-14)\oplus R^{11}(-15)\oplus \rightarrow R^{36}(-13)\oplus R^{11}(-14)\rightarrow
\]
\[
R^{35}(-11)\oplus R^{65}(-12)\oplus R^{10}(-13) \rightarrow R^{10}(-10)\oplus R^{49}(-11) \oplus R^{4}(-12) \rightarrow
\]
\[
R(-7)\oplus R(-8)\oplus R^{35}(-9)\oplus R^{14}(-10)\rightarrow R\rightarrow R/I_{c}(St_{7})\rightarrow 0
\]
is the minimal $\ZZ^{15}$-graded free resolution of $R/I_{c}(St_{7})$. Therefore,
\[
\projdim(R/I_{c}(St_{7}))=\projdim(R/I_{c}(St_{4}))+(\projdim(R/I_{c}(St_{4}))-2)=6.
\]
Then for every $w\geq 1$ we have

\setlength{\tabcolsep}{1.5pt}
\renewcommand{\arraystretch}{1}
\begin{center}
\begin{tabular}{|c|c|}
\hline
$r=3w+1$&$\projdim(R/I_{c}(St_{r}))$\\[1mm]
\hline
4&4\\
\hline
7&$\projdim(R/I_{c}(St_{4}))+(\projdim(R/I_{c}(St_{4}))-2)=4+2=6$\\
\hline
10&$\projdim(R/I_{c}(St_{7}))+(\projdim(R/I_{c}(St_{4}))-2)=6+2=8$\\
\hline
13&$\projdim(R/I_{c}(St_{10}))+(\projdim(R/I_{c}(St_{4}))-2)=8+2=10$\\
\hline
16&$\projdim(R/I_{c}(St_{13}))+(\projdim(R/I_{c}(St_{4}))-2)=10+2=12$\\
\hline
19&$\projdim(R/I_{c}(St_{16}))+(\projdim(R/I_{c}(St_{4}))-2)=12+2=14$\\
\hline
22&$\projdim(R/I_{c}(St_{19}))+(\projdim(R/I_{c}(St_{4}))-2)=14+2=16$\\
\hline
$\vdots$&\vdots\\
\hline
\end{tabular}
\end{center}

\medskip

\noindent Suppose that the induction hypothesis holds for $r-3$. Then
\begin{eqnarray*}
\projdim(R/I_{c}(St_{3w+1}))&=&\projdim(R/I_{c}(St_{3w+1-3}))+ \left(\projdim(R/I_{c}(St_{4}))-2\right)\\
&=&\left(2(w-1)+2\right)+(4-2)\\
&=&2w+2.
\end{eqnarray*}
Then if $r=3w+1$ we have $\projdim(R/I_{c}(St_{r}))=\frac{2}{3}r+\frac{4}{3}$ for every $w \geq1$.
Hence, the assertion follows.
\end{proof}

Finally, we investigate the regularity of $R/I_{c}(St_{r})$.

\begin{Theorem}
\label{reg}
Let $I_{c}(St_{r})$ be the vertex cover ideal of a bipartite planar graph $St_{r}$ with $r > 1$ regions. Then
$$\reg(R/I_{c}(St_{r}))=\left\{ \begin{array}{ll}
\frac{4}{3}r+\frac{1}{3} & \quad \textrm{if \quad  $r=3w-1$},\\[2mm]
\frac{4}{3}r &  \quad \textrm{if \quad $r=3w$},\\[2mm]
\frac{4}{3}r-\frac{1}{3} & \quad \textrm{if \quad $r=3w+1$}.\\
\end{array} \right.
$$
\end{Theorem}

\begin{proof}
Let $St_{r}$ be the bipartite planar graph with $r > 1$ regions, and let
$I_{c}(St_{r})$ be the vertex cover ideal of $St_{r}$.
Hence, the vertex set of $St_{r}$ can be partitioned into disjoint subsets $V_{1}$ and $V_{2}$
such that $N=|V_{1}|+|V_{2}|=2r+1$, and $|V_{i}|$ denotes the number of vertices of $V_{i}$ for $i=1,2$.
Therefore, the minimal  $\ZZ^N$-graded free $R$-resolution of $\reg(R/I_{c}(St_{r}))$ is the following:
\[
 0\to F_p\to F_{p-1} \to \cdots \to F_2\to F_1\to F_0\to R/I_{c}(St_{r})\to 0
\]
where
$
F_p=\Dirsum_{g\in \ZZ}R(-a_{pg}),
$
and
$a_{pg}=(a_{pg}(1),\dots,a_{pg}(N))\in \ZZ^{N}_{+}$. It follows that
$|a_{pg}|=N$ whereas $ a_{pg}(i)=1$ for all $i=1,\dots,N$.

\noindent We have three cases:

\begin{itemize}
\item If $r=3w-1$, then from Theorem \ref{projdim},\,
$
p=\projdim(R/I_{c}(St_{r})=\frac{2}{3}r+\frac{2}{3},
$
\end{itemize}
hence
\begin{eqnarray*}
\reg(R/I_{c}(St_{r}))
&=&N-\projdim(R/I_{c}(St_{r}))\\
&=&(2r+1)-\frac{2}{3}(r+1)\\
&=&\frac{4}{3}r+\frac{1}{3}.
\end{eqnarray*}

\begin{itemize}
\item If $r=3w$, hence Theorem \ref{projdim} implies that
$
p=\frac{2}{3}r+1.
$
It follows that
\end{itemize}
\begin{eqnarray*}
\reg(R/I_{c}(St_{r}))
&=&N-\projdim(R/I_{c}(St_{r}))\\
&=&N-p\\
&=&(2r+1)-\left(\frac{2}{3}r+1\right)\\
&=&\frac{4}{3}r.
\end{eqnarray*}

\begin{itemize}
\item If $r=3w+1$, it follows from Theorem \ref{projdim} that
$
p=\frac{2}{3}r+\frac{4}{3}.
$
Hence,
\end{itemize}
\begin{eqnarray*}
\reg(R/I_{c}(St_{r}))
&=&N-\projdim(R/I_{c}(St_{r}))\\
&=&N-p\\
&=&(2r+1)-\left(\frac{2}{3}r+\frac{4}{3}\right)\\
&=&\frac{4}{3}r-\frac{1}{3}.
\end{eqnarray*}
Then, the assertion follows.
\end{proof}

\begin{Example}
\label{St3}
{\em
Let $St_{3}$ be the bipartite planar graph on vertex set\\[1mm] $V(St_{3})=\{v_1,v_2,v_3,v_4,v_5,v_6,v_7\}$ and
edge set\\[1mm]
$
E(St_{3})=\{\{v_1,v_2\},\{v_1,v_3\},\{v_1,v_4\},\{v_2,v_5\},\{v_3,v_6\},\{v_4,v_7\},\{v_3,v_5\},\{v_4,v_6\},\{v_2,v_7\}\}.\\[-2mm]
$

\begin{center}
\begin{minipage}{.4\textwidth}
\usetikzlibrary {positioning}
\definecolor {processblue}{cmyk}{0.96,0,0,0}
\begin{tikzpicture}
      \tikzset{enclosed/.style={draw, circle, inner sep=0pt, minimum size=.15cm, fill=black, text=blue , minimum width =0.1 cm,node distance =0.1 cm,top color =blue , bottom color = processblue!20 ,
draw,processblue , text=blue , minimum width =0.1 cm }}
\node[enclosed, label={left, yshift=.2cm: $v_{5}$}] (v_{5}) at (2.3,4.5) {};
\node[enclosed, label={left, yshift=.2cm: $v_{2}$}] (v_{2}) at (4,5) {};
\node[enclosed, label={left, yshift=.4cm: $v_{7}$}] (v_{7}) at (5.7,4.5) {};
\node[enclosed, label={left, yshift=.2cm: $v_{3}$}] (v_{3}) at (2,3) {};
\node[enclosed, label={left, yshift=.2cm: $v_{1}$}] (v_{1}) at (4,3.5) {};
\node[enclosed, label={left, yshift=.3cm: $v_{4}$}] (v_{4}) at (6,3) {};
\node[enclosed, label={left, yshift=-.08cm: $v_{6}$}] (v_{6}) at (4,1.5) {};

\draw (v_{2}) -- (v_{5}) node[midway, right] (edge1){};
\draw (v_{2}) -- (v_{7}) node[midway, below] (edge2){};
\draw (v_{2}) -- (v_{1}) node[midway, below] (edge3){};
\draw (v_{5}) -- (v_{3}) node[midway, below] (edge4){};
\draw (v_{1}) -- (v_{3}) node[midway, below] (edge5){};
\draw (v_{1}) -- (v_{4}) node[midway, below] (edge6){};
\draw (v_{4}) -- (v_{7}) node[midway, below] (edge7){};
\draw (v_{3}) -- (v_{6}) node[midway, below] (edge8){};
\draw (v_{4}) -- (v_{6}) node[midway, below] (edge9){};
\end{tikzpicture}
\end{minipage}
\end{center}
\smallskip
The vertex cover ideal $I_{c}(St_{3})$ has $5$ monomial squarefree generators and it is:
\[
I_{c}(St_{3})=(x_{2}x_{3}x_{4},x_{1}x_{5}x_{6}x_{7},x_{1}x_{2}x_{3}x_{6}x_{7},x_{1}x_{3}x_{4}x_{5}x_{7},x_{1}x_{2}x_{4}x_{5}x_{6}).
\]
Then, from Theorems \ref{height}, \ref{projdim} and \ref{reg} it results:
\begin{itemize}
\item[(\.{i})] $\projdim(R/I_{c}(St_{3}))=3$;

\item[(\.{i}\.{i})] $\reg(R/I_{c}(St_{3}))=4$;

\item[(\.{i}\.{i}\.{i})] $\height(I_{c}(St_{3}))^{k}=2$, for every $k \geq1$.
\end{itemize}
}
\end{Example}

\medskip

\section{Syzygy modules of edge ideals of graphs $St_r$}
\label{two}

Let $St_{r}$ be a bipartite planar graph with $r > 1$ regions, and let $I(St_{r})$ be its edge ideal.
In this section we are interested to find the graded Betti numbers that appear in the minimal
graded resolution of $I(St_{r})$.

\begin{Definition}
{\em
Let $G$ be a graph on vertex set $V(G)=\{v_1,\dots,v_n\}$ and let
$E(G)=\{f_1,\dots,f_q\}$ be its edge set. The {\em edge graph} of $G$, denoted by
$\mathcal{L}(G)$, has vertex set equal to the edge set of $G$ and two vertices of $\mathcal{L}(G)$ are
adjacent whenever the corresponding edges of $G$ have one common vertex. It follows that
$V (\mathcal{L}(G))=E(G)=\{f_1,\dots,f_q\}$ and
\begin{equation*}
E(\mathcal{L}(G))=\{(f_i,f_j)\mid f_i=\{v_i,v_j\},  f_j=\{v_j,v_k\}, \quad  i \neq j, j \neq k\}.
\end{equation*}
Then the number of edges of the edge graph $\mathcal{L}(G)$ is given by
\[
|E(\mathcal{L}(G))|=-|E|+ \sum_{i=1}^{n}\frac{\deg^{2}(v_{i})}{2},
\]
where $\deg(v_{i})$ is the number of edges incident with the vertex $v_{i}$ (\cite{Villarreal2015}).
}
\end{Definition}

An interesting problem is to express the second Betti number of $R/I(G)$ in terms of graph theoretical
properties. In \cite{EliahouVillarreal1999}, Eliahou and Villarreal give an explicit formula to compute the
second graded Betti number in degree $3$, that represents the number of the
generators of linear syzygies of $I(G)$.

\begin{Theorem}
\label{second}
Let $I(St_{r})$ be the edge ideal of a bipartite planar graph $St_{r}$ with $r > 1$ regions.
If
\[
\cdots \rightarrow R^{c}(-4)\oplus R^{b}(-3) \rightarrow R^{q}(-2) \rightarrow I(St_{r})\rightarrow 0
\]
is the minimal graded resolution of $I(St_{r})$, it follows that
\begin{itemize}
\item[(\.{i})] $q=3r$;\\[-2mm]

\item[(\.{i}\.{i})] $b=\frac{1}{2}r(r+7)$.
\end{itemize}
\end{Theorem}

\begin{proof}
(\.{i})
Let $St_{r}$ be the planar graph with $r > 1$ regions on vertex set
$V(St_{r})=\{v_1,\dots,v_{2r+1}\}$, and let $I(St_{r})$ be the edge ideal of $St_{r}$.
Therefore,
\begin{eqnarray*}
q=|E|=|\{\{v_1,v_i\} \mid  2 \leq i \leq r + 1\}|+|\{\{v_i,v_{i+r}\} \mid  2 \leq i \leq r + 1\}|
\\
+|\{\{v_i,v_{i+r-1}\} \mid  3\leq i \leq r+1\}|+|\{v_2,v_{2r+1}\}|,
\end{eqnarray*}
hence, $q=r+r+(r-1)+1=3r$.\\

(\.{i}\.{i})
From Eliahou and Villarreal's formula, we have
$
b=|E(\mathcal{L}(St_{r}))|-N_{3},
$
where $N_{3}$ is the number of the triangles of $St_{r}$. Since $St_{r}$ is bipartite, it follows that $N_{3}=0$. Hence,
\begin{eqnarray*}
|E(\mathcal{L}(St_{r}))|&=&-|E(St_{r})|+ \sum_{i=1}^{2r+1}\frac{\deg^{2}v_{i}}{2}
\\
&=&-3r+\frac{r^{2}}{2}+r\left(\frac{3^{2}}{2}\right)+r\left(\frac{2^{2}}{2}\right)
\\
&=&-3r+\frac{1}{2}(r^{2}+13r)
\\
&=&\frac{1}{2}r(r+7),
\end{eqnarray*}
as desired.
\end{proof}

\begin{Example}
\label{St2,syzygy}
{\em
Let $St_{2}$ be the complete bipartite planar graph on vertex set $V(St_{2})=\{v_1,\dots,v_5\}$ as in Example \ref{St2}.
Then
\[
\cdots \rightarrow R^{b}(-3)\rightarrow R^{q}(-2)\rightarrow  I(St_{2})\rightarrow 0.
\]
is the minimal $\ZZ^{5}$-graded free resolution of $I(St_{2})$, where $q=3r=6$ and $b=\frac{1}{2}r(r+7)=9$, according to Theorem \ref{second}.
}
\end{Example}

Furthermore, for any graph $G$, Eliahou and Villarreal in \cite{EliahouVillarreal1999} give a conjecture to compute the third Betti number in degree $4$ of $R/I(G)$, that is the number of the linear generators
of the second syzygy module of $I(G)$.

\medskip

In the following, we prove the conjecture for complete graphs, complete bipartite graphs and cycle graphs.

\begin{Theorem}
\label{complete graph}
Let $\mathcal{K}_{n}$ be a complete graph on vertex set $\{v_1,\dots,v_n\}$, and
$N_{4}$ be the number of the complete subgraphs of $\mathcal{K}_{n}$ on $4$ vertices. If
\[
\cdots \rightarrow \cdots \oplus R^{d}(-4) \rightarrow R^{c}(-4)\oplus R^{b}(-3) \rightarrow R^{q}(-2) \rightarrow I(\mathcal{K}_{n})\rightarrow 0
\]
is the minimal graded free resolution of $I(\mathcal{K}_{n})$, then
\[
d=\sum_{i=1}^{n}\binom{\deg(v_{i})}{3}-N_{4},
\]
where $\deg (v_i)$ is the number of the edges incident with the vertex $v_i$.
\end{Theorem}

\begin{proof}
Let's assume that $\mathcal{K}_{n}$ is a complete graph with vertex set $ V(\mathcal{K}_{n})=\{v_1,\dots,v_n\}$.
The edge ideal of $\mathcal{K}_{n}$ is $
I(\mathcal{K}_{n})=(x_ix_j \mid 1 \leq i < j \leq n).$\,
If $d$ indicates the third graded Betti number $b_{3,4}$ of $R/I(\mathcal{K}_{n})$, it is given
by the formula (\cite{Jacques2004})
\[
d=b_{3,4}(\mathcal{K}_{n})=3\binom{n}{3+1}=3\binom{n}{4}.
\]
On the other hand, $\deg(v_{i})=n-1$ for all $v_{i}\in V(\mathcal{K}_{n})$, $1 \leq i\leq n$.
Hence, $\sum_{i=1}^{n}\binom{\deg(v_{i})}{3}=n\binom{n-1}{3}$ and $N_{4}=\binom{n}{4}$.
Therefore,
\begin{eqnarray*}
\sum_{i=1}^{n}\binom{\deg(v_{i})}{3}-N_{4}&=&n\binom{n-1}{3}-\binom{n}{4}\\
&=&\frac{n!}{3!(n-4)!}-\frac{n!}{4!(n-4)!}\\
&=&3 \left(\frac{n!}{4!(n-4)!}\right)\\
&=&3\binom{n}{4}\\
&=&d,
\end{eqnarray*}
and the assertion follows.
\end{proof}

\begin{Theorem}
\label{complete bipartite}
Let $\mathcal{K}_{n,m}$ be a complete bipartite graph on $N=n+m$
vertices, and $W_{4}$ be the number of squares without chords of $\mathcal{K}_{n,m}$. If
\[
\cdots \rightarrow \cdots\oplus R^{d}(-4) \rightarrow R^{c}(-4)\oplus R^{b}(-3) \rightarrow R^{q}(-2) \rightarrow I(\mathcal{K}_{n,m})\rightarrow 0
\]
is the minimal graded free resolution of $I(\mathcal{K}_{n,m})$, then
\[
d=\sum_{i=1}^{N}\binom{\deg(v_{i})}{3}+W_{4},
\]
where $\deg(v_{i})$ is the number of the edges incident with the vertex $v_{i}$ of $\mathcal{K}_{n,m}$.
\end{Theorem}

\begin{proof}
Let $\mathcal{K}_{n,m}$ be a complete bipartite graph on vertex set $V(\mathcal{K}_{n,m})=V_{1}\union V_{2}=
\{x_1,\dots,x_n\}\union \{y_1,\dots,y_m\}$.
If we assume that $d$ is the third graded Betti number $b_{3,4}$ of $R/I(\mathcal{K}_{n,m})$, it is given by the formula \cite{Jacques2004}
\[
d=b_{3,4}(\mathcal{K}_{n,m})=\binom{n}{1} \binom{m}{3}+\binom{n}{2} \binom{m}{2}+\binom{n}{3} \binom{m}{1}.
\]
Therefore,
\[
\sum_{i=1}^{N}\binom{\deg(v_{i})}{3}=\sum_{j=1}^{n}\binom{\deg(x_{j})}{3}+\sum_{k=1}^{m}\binom{\deg(y_{k})}{3}
=n\binom{m}{3}+m\binom{n}{3}.
\]
Since the squares without chords of $\mathcal{K}_{n,m}$ have two vertices in $V_{1}$ and two vertices in $V_{2}$,
it follows that $W_{4}=\binom{n}{2} \binom{m}{2}$. Hence,
\begin{eqnarray*}
\sum_{i=1}^{N}\binom{\deg(v_{i})}{3}+W_{4}=n\binom{m}{3}+m\binom{n}{3}+\binom{n}{2} \binom{m}{2}=d,
\end{eqnarray*}
and the assertion follows.
\end{proof}

\begin{Theorem}
\label{cycle}
Let $C_{n}$ be a cycle graph with vertex set $\{v_{1},\dots,v_{n}\}$, \,$n > 3$.
If
\[
\cdots \rightarrow \cdots\oplus R^{d}(-4) \rightarrow R^{c}(-4)\oplus R^{b}(-3) \rightarrow R^{q}(-2) \rightarrow I(C_{n})\rightarrow 0
\]
is the minimal graded free resolution of $I(C_{n})$, then
\begin{itemize}
\item[(\.{i})] $d=1\;\;$ for $n=4$;\\[-3mm]

\item[(\.{i}\.{i})] $d=0\;\;$ for $n> 4$.
\end{itemize}
\end{Theorem}

\begin{proof}
(\.{i})
Let $C_{n}$ be a cycle graph with vertex set $V(C_n)=\{v_{1},\dots,v_{n}\}$, and let $I(C_{n})$ be the edge ideal of $C_{n}$.
We consider $n=4$. Then $d=b_{3,4}(C_{4})=1$ (\cite{Jacques2004}).
Since $\deg(v_{i})=2$ for all $v_{i}\in V(C_{4})$, $1\leq i\leq 4$, it follows that
$N_{4}=0$ and $W_{4}=1$. Therefore,
$
\sum_{i=1}^{4}\binom{\deg(v_{i})}{3}-N_{4}+W_{4}=1=d.\\[-2mm]
$

(\.{i}\.{i})
If $n > 4$ we have $d=b_{3,4}(C_{n})=\frac{n}{n-2} \binom{2}{1} \binom{n-2}{1}=0$ (\cite{Jacques2004}).\\
Hence,
$
\sum_{i=1}^{n}\binom{\deg(v_{i})}{3}-N_{4}+W_{4}=0=d,
$
as desired.
\end{proof}
\smallskip
In general, for bipartite graphs, Eliahou and Villarreal's conjecture is stated as follows.\\[2mm]
\textbf{Conjecture.}
{\em
Let $G$ be a bipartite graph on $N$ vertices, $W_{4}$ be the number of squares without chords of $G$.
If
\[
\cdots \rightarrow \cdots\oplus R^{d}(-4) \rightarrow R^{c}(-4)\oplus R^{b}(-3) \rightarrow R^{q}(-2) \rightarrow I(G)\rightarrow 0
\]
is the minimal graded free resolution of the edge ideal $I(G)$, then
\[
d=\sum_{i=1}^{N}\binom{\deg(v_{i})}{3}+W_{4},
\]
where $\deg(v_{i})$ is the number of the edges incident with the vertex $v_{i}$ of $G$.\hfill $\square$
}
\bigskip

Next, we prove the conjecture for bipartite planar graphs $St_{r}$.

\begin{Theorem}
\label{Conjecture}
Let $St_{r}$ be a bipartite planar graph on vertex set $V(St_{r})=\{v_1,\dots,$ $v_{2r+1}\}$, and let $I(St_{r})$ be the edge ideal of $St_{r}$.
Furthermore, we assume that $W_{4}$ be the number of squares without chords of $St_{r}$.
If
\[
\cdots \rightarrow \cdots\oplus R^{d}(-4) \rightarrow R^{c}(-4)\oplus R^{b}(-3) \rightarrow R^{q}(-2) \rightarrow I(St_{r})\rightarrow 0
\]
is the minimal graded free resolution of $I(St_{r})$, then
\[
d=\sum_{i=1}^{2r+1}\binom{\deg(v_{i})}{3}+W_{4},
\]
where $\deg(v_{i})$ is the number of the edges incident with the vertex $v_{i}$.
\end{Theorem}

\begin{proof}
Let $I(St_{r})$ be the edge ideal of $St_{r}$ with $r > 1$ regions.
The vertex set of $St_{r}$ can be partitioned into disjoint subsets $V_{1}$ and $V_{2}$
such that $|V_{1}|+|V_{2}|=2r+1$, and $|V_{i}|$ denotes the number of vertices of $V_{i}$ for $i=1,2$.\\
We prove the assertion by induction on $r$.\\[-3.5mm]
\begin{itemize}
\item If $r=2$, then $St_{2}$ is a complete bipartite planar graph on vertex set
\end{itemize}
$
V(St_{2})=\{v_1,\dots,v_5\}
$
and edge set
\[
E(St_{2})=\{\{v_1,v_2\},\{v_1,v_3\},\{v_2,v_4\},\{v_3,v_5\},\{v_3,v_4\},\{v_2,v_5\}\}.
\]
The minimal $\ZZ^{5}$-graded free resolution of $I(St_{2})$ is the following:
\[
0\rightarrow R(-5)\rightarrow R^{5}(-4)\rightarrow R^{9}(-3)\rightarrow R^{6}(-2)\rightarrow  I(St_{2})\rightarrow 0.
\]
Since $St_{2}$ has also vertex set $V(St_{2})=\{x_1, x_2, x_3\}\union \{y_1, y_2\}$, it follows from \cite{Jacques2004} that
\[
d=b_{3,4}(St_{2})=\binom{3}{1} \binom{2}{3}+\binom{3}{2} \binom{2}{2}+\binom{3}{3} \binom{2}{1}=5.
\]
Hence,
$
\displaystyle{\sum_{i=1}^{5}\binom{\deg(v_{i})}{3}=3\binom{2}{3}+2\binom{3}{3}=2},
$
and $\displaystyle{W_{4}=\binom{3}{2}\binom{2}{2}=3}$. Therefore,
\[
\sum_{i=1}^{5}\binom{\deg(v_{i})}{3}+W_{4}=2+3=5=d.
\]
\vspace{-1mm}
\begin{itemize}
\item If $r=3$, hence $St_{3}$ is the bipartite planar graph on vertex set
\end{itemize}
$
V(St_{3})=\{v_1,\dots,v_{7}\}
$
and edge set
\begin{eqnarray*}
&E(St_{3})=&\{\{v_1,v_i\} \mid  2 \leq i \leq 4\} \cup \{\{v_i,v_{i+3}\} \mid  2 \leq i \leq 4\}
\\
&&
\cup \{\{v_i,v_{i+2}\} \mid  3\leq i \leq 4\}\cup \{v_2,v_{7}\}.
\end{eqnarray*}
Then the minimal $\ZZ^{7}$-graded free resolution of $I(St_{3})$ is the following:
\[
0 \rightarrow R^{2}(-7) \rightarrow R^{8}(-6)\rightarrow R^{9}(-5)\oplus R^{7}(-4) \rightarrow
\]
\[
R^{3}(-4)\oplus R^{15}(-3)\rightarrow R^{9}(-2)\rightarrow I(St_{3})\rightarrow 0.
\]
\noindent
It follows that $d=b_{3,4}(St_{3})=7$. Then
$
\sum_{i=1}^{7}\binom{\deg(v_{i})}{3}=\binom{3}{3}+3\binom{3}{3}+3\binom{2}{3}=4,
$
where $\deg(v_{1})=3$, $\deg(v_{i})=3$ for $2 \leq i \leq 4$ and $\deg(v_{i})=2$ for $5 \leq i \leq 7$.
Therefore, the number of squares without chords of $St_{3}$ is $W_{4}=3$, and hence
\[
\sum_{i=1}^{7}\binom{\deg(v_{i})}{3}+W_{4}=4+3=7=d.
\]
\vspace{-2mm}
\begin{itemize}
\item If $r=4$, hence $St_{4}$ is the bipartite planar graph on vertex set
\end{itemize}
$
V(St_{4})=\{v_1,\dots,v_{9}\}
$
and edge set
\begin{eqnarray*}
&E(St_{4})=&\{\{v_1,v_i\} \mid  2 \leq i \leq 5\} \cup \{\{v_i,v_{i+4}\} \mid  2 \leq i \leq 5\}
\\
&&
\cup \{\{v_i,v_{i+3}\} \mid  3\leq i \leq 5\}\cup \{v_2,v_{9}\}.
\end{eqnarray*}
Then
\[
0 \rightarrow R(-9) \rightarrow R(-8)\oplus R^{12}(-7) \rightarrow R^{36}(-6)\oplus R(-5) \rightarrow
\]
\[
R^{36}(-5)\oplus R^{12}(-4) \rightarrow R^{12}(-4)\oplus R^{22}(-3) \rightarrow R^{12}(-2)\rightarrow I(St_{4})\rightarrow 0
\]
\noindent
is the minimal $\ZZ^{9}$-graded free resolution of $I(St_{4})$. Hence, $d=b_{3,4}(St_{4})=12$.
Therefore,
$
\sum_{i=1}^{9}\binom{\deg(v_{i})}{3}=\binom{4}{3}+4\binom{3}{3}+3\binom{2}{3}=8,
$
where $\deg(v_{1})=4$, $\deg(v_{i})=3$ for $2 \leq i \leq 5$ and $\deg(v_{i})=2$ for $6 \leq i \leq 9$.
Since the number of squares without chords of $St_{4}$ is $W_{4}=4$, it follows that
\[
\sum_{i=1}^{9}\binom{\deg(v_{i})}{3}+W_{4}=8+4=12=d.
\]
Then $d=b_{3,4}(St_{r})$ for $r> 1$ is the following:

\setlength{\tabcolsep}{2.5pt}
\renewcommand{\arraystretch}{1.5}
\begin{center}
\begin{tabular}{|c|c|}
\hline
$r$&$d=b_{3,4}(St_{r})$\\[1mm]
\hline
2&$3\binom{2}{3}+2\binom{3}{3}+(2+1)=5$\\
\hline
3&$\binom{3}{3}+3\binom{3}{3}+3\binom{2}{3}+3=7$\\
\hline
4&$\binom{4}{3}+4\binom{3}{3}+4\binom{2}{3}+4=12$\\
\hline
5&$\binom{5}{3}+5\binom{3}{3}+5\binom{2}{3}+5=20$\\
\hline
6&$\binom{6}{3}+6\binom{3}{3}+6\binom{2}{3}+6=32$\\
\hline
7&$\binom{7}{3}+7\binom{3}{3}+7\binom{2}{3}+7=49$\\
\hline
8&$\binom{8}{3}+8\binom{3}{3}+8\binom{2}{3}+8=72$\\
\hline
$\vdots$&\vdots\\
\hline
\end{tabular}
\end{center}
\vspace{1mm}
Assume that induction hypothesis holds for $r-1$, where $r\geq 4$. It follows that
\[
d=b_{3,4}(St_{r})=b_{3,4}(St_{r-1})+ \binom{r}{3}-\binom{r-1}{3}+2.
\]
Then we have
\[
b_{3,4}(St_{r-1})=\sum_{g=1}^{2r-1}\binom{\deg(v_{g})}{3}+(r-1),
\]
where $\deg(v_{1})=r-1$, $\deg(v_{g})=3$ for $2 \leq g \leq r$ and $\deg(v_{g})=2$ for $r+1 \leq g \leq 2r-1$.
Therefore,
\begin{eqnarray*}
b_{3,4}(St_{r-1})&=&\binom{r-1}{3}+(r-1)\binom{3}{3}+(r-1)\binom{2}{3}+(r-1)\\
&=&\binom{r-1}{3}+(r-1)+(r-1)\\
&=&\binom{r-1}{3}+2(r-1),
\end{eqnarray*}
and hence
\begin{eqnarray*}
d=b_{3,4}(St_{r})&=&\binom{r-1}{3}+2(r-1)+ \binom{r}{3}-\binom{r-1}{3}+2\\
&=&\binom{r}{3}+2r\\
&=&\sum_{i=1}^{2r+1}\binom{\deg(v_{i})}{3}+r,
\end{eqnarray*}
where $\deg(v_{1})=r$; $\deg(v_{i})=3,$ $2 \leq i \leq r+1$; $\deg(v_{i})=2,$ $r+2 \leq i \leq 2r+1$.
\end{proof}

\end{document}